\documentclass[12pt]{amsart}
\advance\hoffset by -0.5 truecm
\usepackage{amssymb}
\newtheorem{Theorem}{Theorem}[section]
\newtheorem{Lemma}[Theorem]{Lemma}
\newtheorem{Corollary}[Theorem]{Corollary}

\newtheorem{Proposition}[Theorem]{Proposition}

\newtheorem{Remark}[Theorem]{Remark}

\def\V{\mbox{Var}}

\def\Div{\mbox{div}}
\def\Z{{\mathbb Z}}
\def\R\re
\def\V{\bf V}

\def \la{\lambda}

\def \re{{\mathbb R}}

\def \K{{\mathbb K}}
\def \V{{\bf V}}

\def \0{\lambda_{0}}
\def \la{\lambda}
\def \ga{\gamma}

\begin{document}
\title[Regularity of weak foliations for thermostats]{Regularity of weak foliations for thermostats}

\author[G.P. Paternain]{Gabriel P. Paternain}
 \address{ Department of Pure Mathematics and Mathematical Statistics,
University of Cambridge,
Cambridge CB3 0WB, England}
 \email {g.p.paternain@dpmms.cam.ac.uk}




\begin{abstract} Let $M$ be a closed oriented surface endowed with a Riemannian metric $g$.
We consider the flow $\phi$ determined by the motion of a particle under the influence of
a magnetic field $\Omega$ and a thermostat with external field ${\bf e}$.
We show that if $\phi$ is Anosov, then it has weak stable and unstable
foliations of class $C^{1,1}$ if and only if the external field ${\bf e}$ has a global potential $U$,
$g_{1}:=e^{-2U}g$ has constant curvature and $e^{-U}\Omega$ is a constant
multiple of the area form of $g_1$.
We also give necessary and sufficient conditions for just one of the weak foliations to be of class
$C^{1,1}$ and we show that
the {\it combined} effect of a thermostat and a magnetic field can produce an Anosov flow with
a weak stable foliation of class $C^{\infty}$ and a weak unstable foliation
which is {\it not} $C^{1,1}$. Finally we study Anosov thermostats depending quadratically on the velocity and we 
characterize those with smooth weak foliations. In particular, we show that quasi-fuchsian flows as defined by Ghys
in \cite{Ghy1} can arise in this fashion.

\end{abstract}

\maketitle

\section{Introduction}
Let $M$ be a closed oriented surface endowed with a Riemannian metric $g$ and let $SM$
be its unit sphere bundle.
Given $f\in C^{\infty}(M,\re)$ and ${\bf e}$ a smooth vector field on $M$ we let $\lambda\in C^{\infty}( SM,\re)$
be the function given by
\[\lambda(x,v):=f(x)+\langle {\bf e}(x),iv\rangle\]
where $i$ indicates rotation by $\pi/2$ according to the orientation
of the surface. The present paper is concerned with the dynamical system determined by the equation
\begin{equation}
\frac{D\dot{\gamma}}{dt}=\la(\gamma,\dot{\ga})\,i\dot{\ga}.
\label{eqgt}
\end{equation}
This equation defines a flow $\phi$ on $SM$ which reduces to the geodesic
flow when ${\bf e}=f=0$.
If $\la$ does not depend on $v$, i.e. ${\bf e}=0$, then $\phi$ is the
{\it magnetic flow} associated with the magnetic field $\Omega:=f\Omega_{a}$, where $\Omega_{a}$ is the area
form of $M$. The magnetic flow can also be seen as the Hamiltonian flow of $|v|_{x}^{2}/2$
with respect to the twisted symplectic form $-d\alpha+\pi^{*}\Omega$, where $\pi:TM\to M$ is the
canonical projection and $\alpha$ is the contact 1-form of the geodesic flow.
If $f=0$, we obtain the {\it Gaussian thermostat} (or isokinetic dynamics, cf. \cite{H}) which is reversible in the
sense that the flip $(x,v)\mapsto (x,-v)$ conjugates $\phi_{t}$ with $\phi_{-t}$ (just as in the case of geodesic flows).
Gaussian thermostats provide interesting models in nonequilibrium statistical mechanics
\cite{Ga,GaRu,Ru1}.

Thus the dynamical system governed by (\ref{eqgt}) describes the motion of a particle on $M$ under the
combined influence of a magnetic field $f\Omega_a$ and a thermostat with external field ${\bf e}$. To avoid cumbersome
terminology we will call $\phi$ simply a {\it thermostat} and we will keep in mind that
a magnetic force may be present. When $f=0$, we will refer to $\phi$ as a Gaussian or
{\it pure} thermostat.

In the present paper we will study the regularity properties of the weak stable and unstable foliations of
{\it Anosov} thermostats. Recall that the Anosov property means that $T(SM)$
splits as $T(SM)=E^{0}\oplus E^{u}\oplus E^{s}$ in such a way that
there are constants $C>0$ and $0<\rho<1<\eta$ such that $E^{0}$ is
spanned by the generating vector field $F$ of the flow,
and for all $t>0$ we have
\[\|d\phi_{-t}|_{E^{u}}\|\leq C\,\eta^{-t}\;\;\;\;\mbox{\rm
and}\;\;\;\|d\phi_{t}|_{E^{s}}\|\leq C\,\rho^{t}.\]
The subbundles are then invariant and H\"older continuous and have
smooth integral manifolds, the stable and unstable manifolds,
which define a continuous foliation with smooth leaves. Also, since we are working
with surfaces, the {\it weak stable} and {\it unstable bundles} given by
\[E^+=\re F\oplus E^s,\]
\[E^-=\re F\oplus E^u\]
are of class $C^{1,\alpha}$, that is, they are $C^1$ and the transverse derivatives
are $\alpha$-H\"older for some $\alpha<1$
depending on the rates of contraction and expansion of the system \cite{HPS,Ha}.
Moreover, if $\phi$ is {\it volume preserving}, the work of S. Hurder and A. Katok \cite{HK}
shows that the weak bundles have transverse derivatives which are Zygmund regular, and thus
$\alpha$-H\"older for any $\alpha<1$.

In general, thermostats are not volume preserving and in \cite{DP2} it is shown that an
Anosov thermostat on a surface preserves an absolutely continuous measure if and only if
${\bf e}$ has global potential, i.e., ${\bf e}=-\nabla U$, where $U\in C^{\infty}(M,\re)$.
If ${\bf e}=-\nabla U$, it is well known (cf. Section 2) that we can make a time change in $\phi$
so that new flow is the magnetic flow of the metric $e^{-2U}g$ and magnetic field
$e^{-U}\Omega$.

The regularity of the strong bundles $E^s$ and $E^u$ for thermostats is now well understood.
Recall that the regularity of the strong bundles could change dramatically
with time changes while the weak bundles do not change at all. Suppose that there exists
a time change so that $E^s\oplus E^u$ is of class $C^1$. Then \cite[Theorem 2.3]{HK}
implies that the reparametrized flow is a contact Anosov flow and in particular volume preserving.
Thus ${\bf e}$ must have a global potential $U$ and the unit sphere bundle of the metric
$e^{-2U}g$ must be a contact type hypersurface in the symplectic manifold
$(TM,-d\alpha+\pi^{*}e^{-U}\Omega)$. The question of which energy levels of a magnetic flow are of contact type
can be decided, at least when $\Omega$ is exact, in terms of Aubry-Mather theory.
For details we refer to \cite{CMP,P}.

If we do not insist on time changes and we just ask when is $E^s\oplus E^u$ of class $C^1$ (or Lipschitz) the situation
is even more rigid and it is described in \cite{P2,DP1,DP2}. The final result is as follows:
if $\Omega$ is exact, then $E^s\oplus E^u$ is of class $C^1$ if and only if ${\bf e}=\Omega=0$ (i.e. $\phi$ is a geodesic
flow); if $\Omega$ is not exact, then $E^s\oplus E^u$ is of class $C^1$ if and only if ${\bf e}=0$, $g$ has constant
curvature and $\Omega$ is a constant multiple of the area form.

Thus the only question that remains to be addressed regarding regularity is that of the {\it weak}
stable and unstable foliations. Here we show:

\medskip

\noindent {\bf Theorem A.} {\it Let $M$ be a closed oriented surface and let
$\phi$ be an Anosov thermostat. Then $\phi$ has weak stable and unstable foliations of class $C^{1,1}$ if and only
if ${\bf e}=-\nabla U$, $g_{1}:=e^{-2U}g$ has constant curvature and $e^{-U}\Omega$ is a constant multiple of the area
form of $g_1$.}

\medskip

Theorem A is based on results of E. Ghys \cite{Ghy1,Ghy2}. Theorem 4.6 in \cite{Ghy2} asserts that
a smooth Anosov flow on a closed 3-manifold with weak stable and unstable foliations of class $C^{1,1}$,
is smoothly orbit equivalent to a suspension or to what Ghys calls a {\it quasi-fuchsian flow}
and which are described in \cite[Th\'eor\`eme B]{Ghy1}. (In our case, since we are working with circles bundles
the latter alternative holds.) A quasi-fuchsian flow $\psi$ depends on a pair of points $([g_1],[g_2])$
in Teichm\"uller space, has smooth weak stable foliation $C^{\infty}$-conjugate to the
weak stable foliation of the constant curvature metric $g_1$ and smooth weak unstable foliation
$C^{\infty}$-conjugate to the weak unstable foliation of the constant curvature metric $g_2$.
Moreover, $\psi$ preserves a volume form if and only if $[g_1]=[g_2]$.

These results of Ghys imply that if a thermostat has weak stable and unstable foliations of class $C^{1,1}$,
then the {\it Godbillon-Vey invariant} of the foliations must be equal to $4\pi^2\chi(M)$.
Following Hurder and Katok \cite{HK} and Y. Mitsumatsu \cite{M} we will compute the Godbillon-Vey
invariants of a thermostat and show that they are equal to $4\pi^2\chi(M)$ if and only if the conditions
in Theorem A hold. What makes the calculation possible is the result proved in \cite{DP2} that asserts
that the weak foliations must be transversal to the vertical fibres of $SM$.

Theorem A is partially motivated by the following surprising result of M. Wojtkowski \cite[Theorem 5.2]{W2}:
if $g$ is a metric of negative curvature on $M$ and ${\bf e}$ a vector field with zero divergence (with respect to $g$),
then the Gaussian thermostat is Anosov, independently of the size of ${\bf e}$ (thus the orbits of an Anosov thermostat
could have very large geodesic curvature). In particular if $g$ has constant negative curvature
and ${\bf e}$ is the vector field dual to a harmonic 1-form, the thermostat is Anosov. Theorem A shows that
such a flow does not have smooth weak foliations.



Since magnetic flows are volume preserving, \cite[Corollary 3.5]{HK} ensures that if one of the weak foliations
is of class $C^{1,1}$, then they are {\it both} in fact $C^{\infty}$. Thus Theorem A combined with this result implies:

\medskip

\noindent {\bf Corollary 1.} {\it Let $M$ be a closed oriented surface and let
$\phi$ be an Anosov magnetic flow. Then $\phi$ has a weak foliation of class $C^{1,1}$ if and only
if $g$ has constant curvature and $\Omega$ is a constant multiple of the area form.}

\medskip

Since Gaussian thermostats are reversible, a similar conclusion can be obtained for them in spite of the fact
that they are not volume preserving:

\medskip

\noindent {\bf Corollary 2.} {\it Let $M$ be a closed oriented surface and let
$\phi$ be an Anosov Gaussian thermostat. Then $\phi$ has a weak foliation of class $C^{1,1}$ if and only
if ${\bf e}=-\nabla U$ and $e^{-2U}g$ has constant curvature.}

\medskip

It is tempting now to speculate that in Theorem A it suffices to assume that just one of the weak foliations
is of class $C^{1,1}$. Surprisingly this is not the case.
As we will see below, the {\it combined} effect of a thermostat and a magnetic field can produce an Anosov flow with
a weak stable foliation of class $C^{\infty}$ and a weak unstable foliation
which is only $C^{1,\alpha}$ with $\alpha<1$ and no more.

In the next result we shall assume for simplicity that ${\bf e}$ is calibrated so that its divergence $\mbox{\rm div} {\bf e}=0$.
As we explain in Remark \ref{zerodiv} this is no restriction at all if one is interested only in the regularity of
the weak foliations.

\medskip

\noindent {\bf Theorem B.} {\it Let $M$ be a closed oriented surface and let
$\phi$ be an Anosov thermostat with $\mbox{\rm div}{\bf e}=0$.
The following are equivalent:
\begin{enumerate}
\item $\phi$ has a weak foliation of class $C^{1,1}$;
\item $\phi$ has a weak foliation of class $C^{\infty}$;
\item $g$ has constant curvature $-c^2$ and there exists $h\in C^{\infty}(M,\re)$
such that $h^2+f^2=c^2$ and
\[{\bf e}=\frac{1}{c^2}i(f\nabla h-h\nabla f).\]
\end{enumerate}
}

\medskip

Of course, if both foliations are of class $C^{1,1}$ we must have ${\bf e}=\nabla f=0$ in agreement with
Theorem A. Theorem B shows that there are thermostats with a $C^{\infty}$ weak stable foliation, but
whose weak unstable foliation is not $C^{1,1}$. Indeed consider a closed hyperbolic surface
$(M,g)$ and let $f$ be a small smooth nonconstant function on $M$. We set
$h:=\sqrt{1-f^2}$ and
\[{\bf e}=\frac{i\nabla f}{\sqrt{1-f^2}}.\]
If $f$ is small enough the thermostat associated with $(g,f,{\bf e})$ will be Anosov and by Theorem B
its weak stable foliation must be of class $C^{\infty}$. On the other hand, its weak unstable
foliation cannot be of class $C^{1,1}$, because otherwise by Theorem A we would have ${\bf e}=\nabla f=0$
which contradicts our choice of $f$.

The results above show that the dissipative quasi-fuchsian flows $\psi$ do not appear as thermostats with
$\la(x,v)=f(x)+\langle {\bf e}(x),iv\rangle$, but if we are prepared to consider a more general $\la$ the situation
changes. Suppose $q$ is a traceless symmetric 2-tensor which we also view as a function on $SM$. Let
\[\la(x,v)=f(x)+\langle {\bf e}(x),iv \rangle+ q_{x}(v,v)\]
and consider the flow $\phi$ on $SM$ defined by (\ref{eqgt}).
Again, to simplify the exposition we will suppose that $\mbox{\rm div}{\bf e}=0$.
Let $\delta_{g}$ be the divergence operator with respect to the metric $g$ acting on symmetric 2-tensors.
Given a symmetric 2-tensor $p$ we can write $p_{x}(u,v)=\langle A_{x}u,v\rangle$ where
$A_{x}:T_{x}M\to T_{x}M$ is a symmetric linear map. Set $\det_{g}p(x):=\det A_{x}$.

\medskip

\noindent {\bf Theorem C.} {\it Let $M$ be a closed oriented surface and let
$\phi$ be an Anosov thermostat defined by
$\la(x,v)=f(x)+\langle {\bf e}(x),iv \rangle+ q_{x}(v,v)$,
where $\mbox{\rm div}{\bf e}=0$.
Then $\phi$ has weak stable and unstable foliations of class $C^{1,1}$ if and only
if
\begin{align*}
&{\bf e}=0,\\
&\delta_{g}(q-fg)=0,\\
&K_{g}+\mbox{\rm det}_{g} (q+fg)=-h^2,
\end{align*}
where $K_g$ is the curvature of $g$ and $h$ is a non-zero constant.
}

\medskip

Let $g_{0}$ be a metric with constant curvature $-1$ and let $q$ be a traceless symmetric 2-tensor
with $\delta_{g_{0}}q=0$. Such a $q$ is the real part of a holomorphic quadratic differential, so they form a real
vector space of dimension $6\,\mbox{\rm genus}(M)-6$. We will see in Section \ref{C} that we can find a unique
metric $g$ in the conformal class of $g_0$ for which
$K_g+\det_{g}(q)=-h^2$ and thus if $q$ is small enough, the thermostat with $\la=q$
is an Anosov flow with $C^\infty$ weak foliations. We will also see that the flow is dissipative
unless $q=0$, so we obtain quasi-fuchsian flows $\psi$ which are not volume preserving.
We do not know if all the quasi-fuchsian flows can be realized in this way, but it is quite likely that
this is the case (see Subsection 6.2).



\section{Preliminaries}

Let $M$ be a closed oriented surface, $SM$ the unit sphere bundle
and $\pi:SM\to M$ the canonical projection. The latter is in fact
a principal $S^{1}$-fibration and we let $V$ be the infinitesimal
generator of the action of $S^1$.

Given a unit vector $v\in T_{x}M$, we will denote by $iv$ the
unique unit vector orthogonal to $v$ such that $\{v,iv\}$ is an
oriented basis of $T_{x}M$. There are two basic 1-forms $\alpha$
and $\beta$ on $SM$ which are defined by the formulas:
\[\alpha_{(x,v)}(\xi):=\langle d_{(x,v)}\pi(\xi),v\rangle;\]
\[\beta_{(x,v)}(\xi):=\langle d_{(x,v)}\pi(\xi),iv\rangle.\]
The form $\alpha$ is the canonical contact form of $SM$ whose Reeb vector
field is the geodesic vector field $X$.

A basic theorem in 2-dimensional Riemannian geometry asserts that
there exists a~unique 1-form $\psi$ on $SM$ (the connection form)
such that $\psi(V)=1$ and
\begin{align}
& d\alpha=\psi\wedge \beta\label{riem1}\\ & d\beta=-\psi\wedge
\alpha\label{riem2}\\ & d\psi=-(K\circ\pi)\,\alpha\wedge\beta
\label{riem3}
\end{align}
where $K$ is the Gaussian curvature of $M$. In fact, the form
$\psi$ is given by
\[\psi_{(x,v)}(\xi)=\left\langle \frac{DZ}{dt}(0),iv\right\rangle,\]
where $Z:(-\varepsilon,\varepsilon)\to SM$ is any curve with
$Z(0)=(x,v)$ and $\dot{Z}(0)=\xi$ and $\frac{DZ}{dt}$ is the
covariant derivative of $Z$ along the curve $\pi\circ Z$.


For later use it is convenient to introduce the vector field $H$
uniquely defined by the conditions $\beta(H)=1$ and
$\alpha(H)=\psi(H)=0$. The vector fields $X,H$ and $V$ are dual to
$\alpha,\beta$ and $\psi$ and as a consequence of (\ref{riem1}--\ref{riem3}) they satisfy the commutation relations
\begin{equation}\label{comm}
[V,X]=H,\quad [V,H]=-X,\quad [X,H]=KV.
\end{equation}
Equations (\ref{riem1}--\ref{riem3}) also imply that the vector fields
$X,H$ and $V$ preserve the volume form $\alpha\wedge d\alpha$ and hence
the Liouville measure.
Note that the flow of $H$ is given by $R^{-1}\circ g_t\circ R$, where $R(x,v)=(x,iv)$
and $g_t$ is the geodesic flow.


Let $\lambda$ be an arbitrary smooth function on $SM$. For several of the results that we will describe below we will not
need $\la$ to be a polynomial of degree $\leq 2$ in the velocities as in the Introduction.
We may still consider a thermostat $\phi$ defined by (\ref{eqgt}) and it is easy to check that
$$
F=X+\lambda V
$$
is the generating vector field of $\phi$.


Now let $\Theta:=-\alpha\wedge d\alpha=\alpha\wedge\beta\wedge\psi$. This volume form generates the
Liouville measure $d\mu$ of $SM$.

\begin{Lemma} We have:
\begin{align}
L_{F}\Theta &=V(\la)\Theta;\label{lie1}\\
L_{H}\Theta &=0;\label{lie2}\\
L_{V}\Theta &=0.\label{lie3}
\end{align}

\label{lied}
\end{Lemma}

\begin{proof} Note that for any vector field $Y$, $L_{Y}\Theta=d(i_{Y}\Theta)$.
Since
$i_{V}\Theta=\alpha\wedge\beta=\pi^*\Omega_{a}$, where $\Omega_{a}$ is the area
form of $M$, we see that $L_{V}\Theta=0$.
Similarly, $L_{X}\Theta=L_{H}\Theta=0$.
Finally $L_{F}\Theta=L_{X}\Theta+L_{\la V}\Theta=d(i_{\la V}\Theta)=
V(\la)\Theta$.

\end{proof}

\subsection{Time changes}

Let $\ga(s)$ be a unit speed solution of
\begin{equation}
\frac{D\gamma'}{ds}=\la(\gamma,\ga')\,i\ga'
\label{termito}
\end{equation}
where $\la(x,v)=f(x)+\langle {\bf e}(x),iv\rangle+q_{x}(v,v)$.
Let $U\in C^{\infty}(M,\re)$ and
\[t(s):=\int_{s_0}^s e^{-U(\ga(\tau))}\,d\tau.\]

\begin{Lemma} Let $\ga_1(t):=\ga(s(t))$ and $g_1:=e^{-2U}g$.
Then $\ga_1(t)$ is a unit speed solution of the thermostat
determined by the quadruple $(g_1,e^{U}f,e^{2U}({\bf e}+\nabla U),e^{-U}q)$.
(Here $\nabla U$ is the gradient of $U$ with respect to the metric $g$.)
\label{repa}
\end{Lemma}

\begin{proof} It is immediate to check that $\ga_1$ has speed one with respect
to $g_1$. Recall that the connection of $g_1$ is given by
\[D^{1}_{X}Y=D_{X}Y-dU(X)Y-dU(Y)X+g(X,Y)\nabla U.\]
Let us indicate derivatives with respect to $t$ with a dot and derivatives with respect to
$s$ with a prime. Then we have
\[\frac{D^{1}\dot{\ga}_{1}}{dt}=\ddot{s}\ga'+\dot{s}^2\frac{D\ga'}{ds}-2\dot{s}^2\langle \nabla U,\ga'\rangle
\ga'+\dot{s}^2\nabla U.\]
Using that $\dot{s}=e^{U}$ and (\ref{termito}) we obtain:
\[\frac{D^{1}\dot{\ga}_{1}}{dt}=e^{2U}(fi\ga'+\langle {\bf e}+\nabla U,i\ga'\rangle i\ga'
+q(i\ga',i\ga')i\ga')\]
and the lemma follows.

\end{proof}

\begin{Remark}{\rm
Note that we may always choose $U$ so that $\mbox{\rm div}({\bf e}+\nabla U)=0$. Since the divergence
of $e^{2U}({\bf e}+\nabla U)$ with respect to the metric $g_1$ is also zero, the lemma ensures
that given a thermostat determined by $(g,f,{\bf e},q)$ we may always perform a smooth time change so that
the reparametrized flow is the flow of a thermostat whose external field has zero divergence.}
\label{zerodiv}
\end{Remark}

\subsection{Properties of Anosov thermostats} In this subsection we collect some
properties of Anosov thermostats $\phi$ for $\la$ arbitrary.
A result of E. Ghys \cite{Ghy} ensures that $\phi$ is topologically
conjugate to the geodesic flow of a metric of constant negative curvature and thus $\phi$ is transitive and
topologically mixing. This fact is used in the proof of Lemma \ref{trans} below.

Recall the definition of the weak stable and unstable bundles:
\[E^+=\re F\oplus E^s,\]
\[E^-=\re F\oplus E^u.\]

\begin{Lemma}\cite{DP2} For any $(x,v)\in SM$, $V(x,v)\notin E^{\pm}(x,v)$.
\label{trans}
\end{Lemma}

The lemma implies that there exist unique functions $r^{\pm}$ on $SM$
such that
\[H+r^+ V\in E^+,\]
\[H+r^- V\in E^-.\]
Note that the Anosov property implies that $r^+\neq r^-$ everywhere.
Below we will need to use that the functions $r^{\pm}$ satisfy a Riccati
type equation along the flow. Note that $r^{\pm}$ are as smooth as
the bundles $E^\pm$ (i.e. $C^{1,\alpha}$).

\begin{Lemma} \cite{DP2} Let $r=r^{\pm}$. Then
\[F(r-V(\la))+r(r-V(\la))+\K=0,\]
where ${\mathbb K}:=K-H(\la)+\la^2+F(V(\la))$.
\label{riccati}
\end{Lemma}

Later on we will need the following integrated version of this equation:

\begin{Lemma} We have:
\[\int_{SM}\{(r-V(\la))^2+\la^2\}\,d\mu=-4\pi^2\chi(M)+\int_{SM}[V(\la)]^2\,d\mu.\]
\label{riccint}
\end{Lemma}

\begin{proof} Since $L_{F}\Theta=V(\la)\Theta$ an easy consequence of the
Stokes' theorem shows that
\[\int_{SM}F(r)\,d\mu=-\int_{SM}rV(\la)\,d\mu\]
hence integrating the equation in Lemma \ref{riccati} we obtain
\[-\int_{SM}rV(\la)\,d\mu+\int_{SM}(r-V(\la))r\,d\mu=
-\int_{SM}(K+\la^2-H(\la))\,d\mu.\]
Since $H$ preserves the Liouville measure
\[\int_{SM}H(\la)\,d\mu=0,\]
and by the Gauss-Bonnet theorem
\[\int_{SM}K\,d\mu=4\pi^{2}\chi(M),\]
so the lemma follows.

\end{proof}

We conclude this subsection with the following simple lemma which follows right away
from Lemma \ref{riccati}. Note that
$r^{+}-r^{-}$ is function of class $C^{1,\alpha}$ which never vanishes and without loss of generality we may assume
it is always positive.

\begin{Lemma}
$$F(\log(r^{+}-r^{-}))=V(\la)-(r^{+}+r^{-}).$$
\label{vol}
\end{Lemma}

\section{The Godbillon-Vey invariant}

We briefly recall the definition of the Godbillon-Vey class and the
Godbillon-Vey invariant. Let $X$ be a smooth manifold and ${\mathcal F}$ a codimension-one foliation
with $C^{\infty}$ leaves and transversally $C^2$. We suppose that the normal bundle to
${\mathcal F}$ is oriented so that there is a 1-form $\tau$ whose kernel coincides with
$T{\mathcal F}$. By the Frobenius theorem there exists a 1-form $\eta$ of class $C^1$
such that $d\tau=\eta\wedge\tau$. The continuous 3-form $\eta\wedge d\eta$ is closed and its cohomology class
$GV({\mathcal F})\in H^3(X,\re)$, called the {\it Godbillon-Vey class},
is an invariant of diffeomorphism and foliated concordance.
When $X$ is a closed oriented 3-manifold, one can also define the {\it Godbillon-Vey invariant}
\[gv({\mathcal F}):=\int_{X}\eta\wedge d\eta.\]

Hurder and Katok \cite{HK} proved that if $\mathcal F$ is a foliation which is transversally only
$C^{1,\alpha}$ for $\alpha>1/2$, then there is a natural, well-defined Godbillon-Vey invariant
$gv({\mathcal F})$, which extends the previous definition for $\mathcal F$ transversally $C^2$.

The next proposition is the natural generalization to thermostats of the Mitsumatsu formula
for geodesic flows, see \cite[Proposition 9.1]{HK} and \cite{M}.

\begin{Proposition} Let $M$ be a closed oriented surface and let $\phi$ be an Anosov thermostat with
$\lambda$ arbitrary.
Let $\mathcal F$ be one of the weak foliations and suppose it is of class $C^{1,\alpha}$ with
$\alpha>1/2$. Then
\[gv({\mathcal F})=4\pi^{2}\chi(M)-3\int_{SM}([V(\la)]^2+[V(r)]^2)\,d\mu
+2\int_{SM}V(r)(V^2(\la)-2\la)\,d\mu\]
where $r$ is the unique function of class $C^{1,\alpha}$ such that $H+rV\in T{\mathcal F}$.
\label{godveyg}
\end{Proposition}

\begin{Remark}{\rm Note that if $\phi$ is a magnetic flow, it is volume preserving and hence
by \cite[Theorem 3.1]{HK} the weak foliations of $\phi$ are of class $C^{1,\alpha}$ for
{\it any} $\alpha<1$ and hence $gv({\mathcal F})$ is well defined. The Godbillon-Vey invariant is also well defined
under an appropriate pinching condition. Let $\eta_s$ be a unit vector field spanning $E^s$ and let $\eta_u$ be a
unit vector field spanning
$E^u$. Define $\la_s$ and $\la_u$ by the equations:
\begin{align*}
&\la_{s}(x,v,t):=\log|d\phi_{t}(\eta_{s}(x,v))|,\\
&\la_{u}(x,v,t):=-\log|d\phi_{t}(\eta_{u}(x,v))|.
\end{align*}
We say that $\phi$ is $\tau$-{\it pinched} if for all $t$ and $(x,v)$ we have
\[\frac{1}{\tau}\la_s(x,v,t)\leq\la_{u}(x,v,t)\leq \tau \la_{s}(x,v,t).\]
This condition does not depend on time changes (cf. \cite[Section 3.3]{Ghy1}).
Note also that $\phi$ is 1-pinched if and only if $\phi$ preserves a volume form.

If $\phi$ is $\tau$-pinched with $\tau<2$, then Corollary 1.8 in \cite{Ha} ensures that the weak
foliations are $C^{1,\alpha}$ with $\alpha>1/2$. Thus any thermostat which is $C^1$-close to a
volume preserving flow will have a well defined Godbillon-Vey invariant. }
\end{Remark}

\begin{proof} We shall calculate $gv({\mathcal F})$ as if $\mathcal F$ were transversally of class
$C^2$. This is really no restriction as we can always approximate the 1-form $\eta$ by a sequence
of $C^1$ forms $\eta_n$ which converge to $\eta$ in the $C^{\alpha}$-topology as in the proof of
Proposition 9.1 in \cite{HK}. In fact, the definition of the extension of the Godbillon-Vey invariant
in \cite{HK} to $C^{1,\alpha}$ foliations with $\alpha>1/2$ is so that the required continuity holds.

By Lemma \ref{trans}, $V$ is transversal to $\mathcal F$ and hence the 1-form
$$\tau:=-\la\alpha-r\beta+\psi$$
vanishes on $T{\mathcal F}$ and takes the value 1 on $V$.
Hence we may take $\eta$ to be the 1-form given by
\[\eta:=-i_{V}d\tau.\]
Using the identities (\ref{riem1}--\ref{riem3}) we calculate:
\[\eta=-(r-V(\la))\alpha+(\la+V(r))\beta.\]
Finally we compute
\[\eta\wedge d\eta=
\{-(\la+V(r))^2-(r-V(\la))^2+(r-V(\la))V(\la+V(r))\]
\[-(\la+V(r))V(r-V(\la))\}
\alpha\wedge\beta\wedge\psi\]
thus
\[gv({\mathcal F})=\int_{SM}\{-(\la+V(r))^2-(r-V(\la))^2+(r-V(\la))V(\la+V(r))\]
\[-(\la+V(r))V(r-V(\la))\}\,d\mu.\]
Using the fact that $V$ preserves $\mu$ we can rewrite the last integral as
\[gv({\mathcal F})=\int_{SM}\{-(\la+V(r))^2-(r-V(\la))^2-2(\la+V(r))V(r-V(\la))\}\,d\mu.\]
Expanding and simplifying we get
\[gv({\mathcal F})=\int_{SM}\{-\la^2-(r-V(\la))^2-3[V(r)]^2
-4\la V(r)+2V^2(\la)(\la+V(r))\}\,d\mu.\]
Using again that $V$ preserves $\mu$ we have
\[\int_{SM}\la V^2(\la)\,d\mu=-\int_{SM}[V(\la)]^2\,d\mu\]
which implies
\[gv({\mathcal F})=\int_{SM}\{-\la^2-(r-V(\la))^2-3[V(r)]^2
+2V(r)(V^2(\la)-2\la)-2[V(\la)]^2\}\,d\mu\]
and the proposition follows from Lemma \ref{riccint}.

\end{proof}

We will now rewrite the formula in the proposition using a bit of Fourier analysis.
Let $L^2(SM)$ be the space of square integrable functions with respect
to the Liouville measure of $SM$.
The space $L^{2}(SM)$ decomposes into an orthogonal direct sum of
subspaces $\sum H_{n}$, $n\in\Z$, such that on $H_{n}$, $-i\,V$ is
$n$ times the identity operator (cf. \cite{GK}). Given $\la\in C^{\infty}(SM)$ we can expand it as follows:
\[\la=\sum_{k=0}^{\infty} Q_{k}\]
where $Q_k$ are smooth functions such that $Q_{k}\in H_{-k}\oplus H_{k}$.
Set $P_{0}:=Q_{0}$ and for $k\geq 1$ set
\[P_{k}:=\frac{-(k^2+2)}{3k^2}V(Q_{k}).\]

\begin{Proposition} Under the same hypothesis as in Proposition \ref{godveyg} we have
\[gv({\mathcal F})=4\pi^2\chi(M)-3\int_{SM}\left[V\left(r+\sum_{k=0}^{\infty}P_k\right)\right]^2\,d\mu\]
\[+3\sum_{k=3}^{\infty}\left[\left(\frac{k^{2}+2}{3}\right)^2-k^2\right]\int_{SM} Q_{k}^2\,d\mu.\]
\label{godveyg2}
\end{Proposition}

\begin{proof} In terms of the $L^2$-inner product, Proposition \ref{godveyg} gives
\[gv({\mathcal F})=4\pi^2\chi(M)-3\left(||V(\la)||^2+||V(r)||^2\right)+2\langle V(r),V^2(\la)-2\la\rangle.\]
Now observe the identities
\[V^2(Q_{k})=-k^2 Q_{k},\]
\[||V(Q_{k})||^2=k^2||Q_{k}||^2,\]
\[V^{2}(\la)-2\la=-\sum_{k=0}^{\infty}(2+k^2)Q_{k}.\]
Using that elements in $H_{-k}\oplus H_{k}$ are orthogonal to elements in $H_{-l}\oplus H_{l}$
for $k\neq l$ it is now straightforward to check that the identity claimed in the proposition holds.

\end{proof}

We now specialize Proposition \ref{godveyg2} to the case that we are interested in, i.e., when
$\lambda=f+V(\theta)+q$, where $\theta$ is the 1-form dual to the external field $E$
and $q\in H_{-2}\oplus H_{2}$.

\begin{Corollary} Let $M$ be a closed oriented surface and let $\phi$ be an Anosov thermostat
defined by $\lambda=f+V(\theta)+q$.
Let $\mathcal F$ be one of the weak foliations and suppose it is of class $C^{1,\alpha}$ with
$\alpha>1/2$. Then
\[gv({\mathcal F})=4\pi^{2}\chi(M)-3\int_{SM}[V(r+\theta-V(q)/2)]^2\,d\mu\]
where $r$ is the unique function of class $C^{1,\alpha}$ such that $H+rV\in T{\mathcal F}$.
\label{godveyt}
\end{Corollary}

\begin{proof} It suffices to note that since $Q_{1}=V(\theta)$ and $Q_{2}=q$, then
\[P_{1}=-V(V(\theta))=\theta\]
and
\[P_{2}=-V(q)/2.\]


\end{proof}

\section{Proof of Theorem A}

Suppose first that ${\bf e}=-\nabla U$, $g_{1}:=e^{-2U}g$ has constant curvature
and $e^{-U}\Omega$ is a constant multiple of the area form of $g_1$.
On account of Lemma \ref{repa} we can make a smooth time change to the flow $\phi$
so that the reparametrized flow is the magnetic flow $\psi$ of the constant curvature metric
$g_1$ and magnetic field given by a constant multiple of the area form of $g_1$.
Hence the flow $\psi$ is algebraic and thus it has smooth weak foliations.
Consequently, $\phi$ has also smooth weak foliations.

Suppose now that $\phi$ has both weak foliations of class $C^{1,1}$. By \cite[Theorem 4.6]{Ghy2}
we can perform a smooth time change on $\phi$ so that the new flow is smoothly conjugate to a
quasi-fuchsian flow as defined in \cite{Ghy1}. Since the Godbillon-Vey invariant of the weak foliations
of a quasi-fuchsian flow equals $4\pi^2\chi(M)$ and the Godbillon-Vey invariant does not change under conjugacies
we conclude that for our thermostat $\phi$ we must have
\[gv^{\pm}=4\pi^2\chi(M).\]
By Corollary \ref{godveyt} this implies that
\[V(r^{\pm}+\theta)=0.\]
If we let $h^{\pm}=r^{\pm}+\theta$, we can think of $h^{\pm}$ as smooth functions
defined on $M$.
By Lemma \ref{vol} we have
\[\theta=F(\log(h^{+}-h^{-}))+h^{+}+h^{-}.\]
But if $g$ is a smooth function on $M$, then $F(g\circ\pi)=dg$ and the last equality
implies right away that $h^{+}+h^{-}=0$ and
\[\theta=F(\log(h^{+}-h^{-})).\]
Thus $\theta$ is an exact 1-form (i.e. ${\bf e}=-\nabla U$ for some smooth function $U$) and by
Lemma \ref{repa} we may perform a time change so that
we just have to deal with a magnetic flow whose functions $r^{\pm}$ satisfy
$V(r^{\pm})=0$. Again, the equality in Lemma \ref{vol} tells us that
$r^{+}+r^{-}=0$ and $F(\log(2r^{+}))=0$ and thus $r^{\pm}$ are constant functions.
If we now input this information into the Riccati type equation from Lemma \ref{riccati}
we obtain
\[H(f\circ\pi)=r^2+K+f^2\]
which implies that $f$ and $K$ are constant functions as desired.

\qed

\section{Proof of Theorem B}

We first characterize the thermostats for which the Godbillon-Vey invariant is maximal.

\begin{Proposition}Let $M$ be a closed oriented surface and let
$\phi$ be an Anosov thermostat with $\mbox{\rm div}E=0$.
Suppose the weak foliations are of class $C^{1,\alpha}$ with
$\alpha>1/2$.
Then there is a weak foliation with $gv=4\pi^2\chi(M)$ if and only if
$g$ has constant curvature $-c^2$ and there exists $h\in C^{\infty}(M,\re)$
such that $h^2+f^2=c^2$ and
\[{\bf e}=\frac{1}{c^2}i(f\nabla h-h\nabla f).\]

\label{gve}
\end{Proposition}

To prove this proposition we need some preparations. Following V. Guillemin and
D. Kazhdan in \cite{GK} we introduce the following first order differential operators:

\[\eta_{+}:=(X-i\,H)/2\]
and
\[\eta_{-}:=(X+i\,H)/2.\]

Let $L^2(SM)$ be the space of square integrable functions with respect
to the Liouville measure of $SM$.
We now summarize some of the main properties of these
operators (cf. \cite{GK}).

\begin{enumerate}
\item $L^{2}(SM)$ decomposes into an orthogonal direct sum of
subspaces $\sum H_{n}$, $n\in\Z$, such that on $H_{n}$, $-i\,V$ is
$n$ times the identity operator;
\item $\eta_{+}$ extends to a densely defined operator from
$H_{n}$ to $H_{n+1}$ for all $n$. Moreover, its transpose is
$-\eta_{-}$;
\item let $C_{n}^{\infty}(SM)=H_{n}\cap C^{\infty}(SM)$. The operators
$\eta_{\pm}:C^{\infty}_{n}\to C^{\infty}_{n\pm 1}$ are first order elliptic
differential operators.
\end{enumerate}

Given a smooth 1-form $\theta$ we can decompose $\theta$ as
\[\theta=\theta_{-1}+\theta_{1}\]
where
\[2\theta_{-1}=\theta+i V(\theta),\]
\[2\theta_{1}=\theta-i V(\theta).\]
Clearly $\theta_{\pm 1}\in H_{\pm 1}$. Note that $\overline{\eta_{+}\theta_{-1}}=\eta_{-}\theta_{1}$.

The following easy lemma will be useful later on.

\begin{Lemma} The form $\theta$ is closed if and only if $\Im \eta_{-}\theta_{1}=0$.
The form $\theta$ is coclosed if and only if $\Re \eta_{-}\theta_{1}=0$.
Also, $\theta$ is closed if and only $V(\theta)$ coclosed (and hence $\theta$ is coclosed
if and only if $V(\theta)$ is closed since $V^2(\theta)=-\theta$).
\label{util}
\end{Lemma}

\begin{proof}

Note that $4\Re \eta_{-}\theta_{1}=X(\theta)+HV(\theta)$ and that
$4\Im \eta_{-}\theta_{1}=H(\theta)-XV(\theta)$.
Let ${\bf e}$ be the vector field dual to $\theta$. Then $\theta_{x}(v)=\langle {\bf e}(x),v\rangle$
and $V(\theta)=\langle {\bf e}(x),iv\rangle$. We now compute
\[X(\theta)(x,v)=\langle \nabla_{v}{\bf e},v\rangle.\]
Using the expression for the flow of $H$ given in Section 2 we also get
\[H(\theta)(x,v)=\langle \nabla_{iv}{\bf e},v\rangle.\]
Thus
\[X(\theta)+HV(\theta)=\langle \nabla_{v}{\bf e},v\rangle+\langle \nabla_{iv}{\bf e},iv\rangle
=\mbox{\rm div}{\bf e}\]
and
\[H(\theta)-XV(\theta)=\langle \nabla_{iv}{\bf e},v\rangle-\langle \nabla_{v}{\bf e},iv\rangle
=-d\theta(v,iv).\]

\end{proof}

\medskip

\noindent{\it Proof of Proposition \ref{gve}.} Suppose first there is a weak foliation
with $gv=4\pi^2\chi(M)$. By Corollary \ref{godveyt} there exists a function
$h\in C^{1}(M,\re)$ such that $r=h\circ\pi-\theta$ is a solution of the Riccati type
equation of Lemma \ref{riccati} (in what follows we write $h$ also for $h\circ\pi$ and similarly for
$f$ and $K$).
Thus
\[F(h)+h(h-\theta)+K-H(f+V(\theta))+(f+V(\theta))^2-F(\theta)=0.\]
Using that $F=X+(f+V(\theta))V$ and that $V(h)=0$ we obtain
\[X(h)-h\theta-H(f)+fV(\theta)+K+h^2+f^2-(X(\theta)+HV(\theta))=0\]
and since $\mbox{\rm div}{\bf e}=0$, $ 4\Re \eta_{-}\theta_{1}=X(\theta)+HV(\theta)=0$ and therefore
\begin{equation}
X(h)-h\theta-H(f)+fV(\theta)+K+h^2+f^2=0.
\label{simplr}
\end{equation}
Now observe that $X(h)-h\theta-H(f)-fV(\theta)\in H_{-1}\oplus H_{1}$ and
$K+h^2+f^2\in H_{0}$, so equation (\ref{simplr}) is equivalent to the pair of equations
\begin{align}
& X(h)-h\theta-H(f)+fV(\theta)=0,\label{imp}\\
& K+h^2+f^2=0.\label{curv}
\end{align}
Consider the complex valued function $p:=h+if$ and let
$\omega:= X(h)-h\theta-H(f)+fV(\theta)$.
A simple calculation shows that
\[\eta_{-}(p)-\theta_{-1}p=(\omega
+iV(\omega))/2=\omega_{-1}.\]
It follows that equation (\ref{imp}) which is just $\omega=0$, is equivalent to
\begin{equation}
\eta_{-}(p)-\theta_{-1}p=0.
\label{cf}
\end{equation}
Since $\eta_{\pm}$ are elliptic operators, equation (\ref{cf}) tells us that
$h\in C^{\infty}(M,\re)$.

\begin{Lemma} If $p$ is a solution of (\ref{cf}) then
\[\eta_{+}\eta_{-}(|p|^2)=|\eta_{+}(p)+\theta_{1}p|^2+2|p|^2\Re \eta_{-}\theta_{1}.\]
\end{Lemma}

\begin{proof} This is an easy calculation in which one uses that
$\overline{\eta_{-}p}=\eta_{+}\overline{p}$ and that $\eta_{+}\eta_{-}=\eta_{-}\eta_{+}$
on $H_0$. We omit the details.
\end{proof}

Since we are assuming that $\mbox{\rm div}{\bf e}=0$, Lemma \ref{util} tells us that
$\Re \eta_{-}\theta_{1}=0$ and therefore, on account of the last lemma
\[\eta_{+}\eta_{-}(|p|^2)=|\eta_{+}(p)+\theta_{1}p|^2\geq 0.\]
Since $\eta_{+}\eta_{-}$ on $H_{0}$ is nothing but the Laplace-Beltrami operator in disguise, the last
inequality implies that $|p|^2$ must be a constant function on $M$
and
\begin{equation}
\eta_{+}(p)+\theta_{1}p=0.
\label{cfc}
\end{equation}
The fact that $|p|^2$ is constant, combined with (\ref{curv}) shows that the metric $g$ must have constant
curvature, let us say, $-c^2$.
If we now combine equations (\ref{cf}) and (\ref{cfc}) we obtain:
\[X(p)-ipV(\theta)=0\]
which is equivalent to
\[X(h)+fV(\theta)=0\]
and
\[X(f)-hV(\theta)=0.\]
Solving these equations for $\theta$ we arrive at
\[c^2\theta=fH(h)-hH(f)\]
which shows that the conditions in the proposition are necessary to have a maximal
Godbillon-Vey invariant.

To show that the conditions are also sufficient we only need to observe that
if the metric has constant curvature $-c^2$ and $h$ is a smooth function
such that $h^2+f^2=c^2$ and
\[c^2\theta=fH(h)-hH(f)\]
then $h$ and $f$ satisfy the pair of equations (\ref{imp}) and (\ref{curv})
and thus if we let $r=h-\theta$, then $r$ is a $C^{\infty}$ function which satisfies
the Riccati type equation of Lemma \ref{riccati}. It follows that
the smooth vector fields $F$ and $H+rV$ span a two dimensional bundle invariant under
the Anosov thermostat. This bundle can only be one of the weak bundles which shows
that there is a $C^{\infty}$ weak foliation which has $gv=4\pi^2\chi(M)$ by Corollary \ref{godveyt}.

\qed

\subsection{Proof of Theorem B} Let us show that (1) implies (3) and let $\mathcal F$ be a weak foliation
of class $C^{1,1}$.
By \cite[Theorem 4.1]{Ghy2}, $\mathcal F$ is transversally projective and thus $\mathcal F$ is
$C^{1}$-conjugate to the weak foliation of a geodesic flow of constant negative curvature (see pages 178 and 179
in \cite{Ghy2}). As in the proof
of Theorem A we can conclude that $gv({\mathcal F})=4\pi^2\chi(M)$ and
Proposition (\ref{gve}) implies that (3) holds.
On the other hand (3) implies (2): it suffices to check that the argument given at the end
of the proof of Proposition (\ref{gve}) shows that if (3) holds then there is a weak bundle which
is $C^{\infty}$. Since (2) obviously implies (1), the proof of Theorem B is now complete.

\section{Proof of Theorem C and complements}\label{C}

\subsection{Proof of Theorem C}

Suppose that $\phi$ has both weak foliations of class $C^{1,1}$. We now argue as in the proof of Theorem A and note
that the presence of $q$ does not really affect the argument that shows that $\theta$ is exact.
By \cite[Theorem 4.6]{Ghy2}
we can perform a smooth time change on $\phi$ so that the new flow is smoothly conjugate to a
quasi-fuchsian flow as defined in \cite{Ghy1}. Since the Godbillon-Vey invariant of the weak foliations
of a quasi-fuchsian flow equals $4\pi^2\chi(M)$ and the Godbillon-Vey invariant does not change under conjugacies
we conclude that for our thermostat $\phi$ we must have
\[gv^{\pm}=4\pi^2\chi(M).\]
By Corollary \ref{godveyt} this implies that
\[V(r^{\pm}+\theta-V(q)/2)=0.\]
If we let $h^{\pm}=r^{\pm}+\theta-V(q)/2$, we can think of $h^{\pm}$ as smooth functions
defined on $M$. Note that $V(\la)=-\theta+V(q)$ and hence
by Lemma \ref{vol} we have
\[\theta=F(\log(h^{+}-h^{-}))+h^{+}+h^{-}.\]
The last equality implies right away that $h^{+}+h^{-}=0$ and
\[\theta=F(\log(h^{+}-h^{-})).\]
Thus $\theta$ is an exact 1-form and since we are assuming that it is coclosed, it must vanish identically.
It follows that $h^{\pm}$ must be constant functions.

Let us now use that $r=h+V(q)/2$, where $h=h^{\pm}$ is a constant, is a solution of the
Riccati type equation of Lemma \ref{riccati} with $\la=f+q$.
We obtain:
\[F(h-V(q)/2)+(h+V(q)/2)(h-V(q)/2)+K-H(f+q)+(f+q)^2+F(V(q))=0\]
where $F=X+(f+q)V$. Using that $h$ is a constant and $V^2(q)=-4q$ we derive:
\begin{equation}
XV(q)/2-H(f+q)+K+h^2+f^2-q^2-[V(q)]^2/4=0.
\label{simplr2}
\end{equation}
Now observe that $XV(q)/2-H(f+q)\in H_{-1}\oplus H_{1}$ and
$K+h^2+f^2-q^2-[V(q)]^2/4\in H_{0}$, so equation (\ref{simplr2}) is equivalent
to the pair of equations
\begin{align}
& XV(q)/2-H(f+q)=0,\label{imp2}\\
& K+h^2+f^2-q^2-[V(q)]^2/4=0.\label{curv2}
\end{align}

It is straightforward to check that
\[\mbox{\rm det}_{g}(q+fg)=f^2-q^2-[V(q)]^2/4\]
and thus equation (\ref{curv2}) can be rewritten as
\[K+h^2+\mbox{\rm det}_{g}(q+fg)=0.\]
To understand equation (\ref{imp2}) we will
use the following easy lemma.

\begin{Lemma} Let $q$ be a symmetric 2-tensor which we also view as a function on $SM$.
Then $\delta_{g}q\,|_{SM}=X(q)+HV(q)/2$.
\label{ldiv}
\end{Lemma}

\begin{proof} By definition of $\delta_g$:
\[\delta_{g}q(x,v)=(\nabla_{v}q)(v,v)+(\nabla_{iv}q)(iv,v).\]
Let $\ga$ be the geodesic with initial conditions $(x,v)$. Then
\[X(q)(x,v)=\left.\frac{d}{dt}\right|_{t=0}q_{\ga(t)}(\dot{\ga}(t),\dot{\ga}(t))
=(\nabla_{v}q)(v,v).\]
Note that
\[V(q)(x,v)=2\,q_{x}(iv,v)\]
hence using the expression for the flow of $H$ given in Section 2 we also get
\[HV(q)(x,v)=2(\nabla_{iv}q)(iv,v).\]

\end{proof}

By applying $V$ we see that equation (\ref{imp2}) is equivalent to
\[X(q-f)+HV(q)/2=0\]
and thus the lemma implies that (\ref{imp2}) is equivalent to
\[\delta_{g}(q-fg)=0.\]
Summarizing, we have shown that (\ref{imp2}) and (\ref{curv2}) are equivalent to
\begin{align}
&\delta_{g}(q-fg)=0, \label{div2}\\
&K_{g}+h^2+\mbox{\rm det}_{g} (q+fg)=0. \label{curvu}
\end{align}
This shows that the equations listed in Theorem C are necessary for having
weak foliations of class $C^{1,1}$. To show that they are sufficient
we note that if we let $r^{\pm}=\pm h+V(q)/2$, then $r^{\pm}$ is a $C^{\infty}$ function
which satisfies the Riccati type equation of Lemma \ref{riccati}. It follows that
the smooth vector fields $F$ and $H+r^{\pm}V$ span a two dimensional bundle invariant under
the Anosov thermostat. This bundle can only be one of the weak bundles which shows
that the weak foliations are $C^{\infty}$ (with $gv=4\pi^2\chi(M)$ by Corollary \ref{godveyt}).
This finishes the proof of Theorem C.

\qed

\subsection{The space of solutions when $f=0$}
In this subsection we consider the solutions of the equations
\begin{align}
&\delta_{g}(q)=0, \label{div2a}\\
&K_{g}+h^2+\mbox{\rm det}_{g} (q)=0, \label{curvua}
\end{align}
which arise from Theorem C when ${\bf e}=f=0$.

Given a $C^{\infty}$ metric $g$ on $M$, its conformal class contains a unique
metric $g_0$ with constant curvature $-1$. The metric $g_0$ determines a complex
structure and the equation $\delta_{g}(q)=0$ simply says that $q$ is the real part of a
holomorphic quadratic differential. Write $g=e^{2u}\,g_0$ and note that
\[\mbox{\rm det}_{g} (q)=e^{-4u}\,\mbox{\rm det}_{g_{0}} (q).\]
The curvature $K_g$ can be expressed in terms of $u$ as
\[K_g=e^{-2u}(\Delta_{g_0}u-1),\]
where $\Delta_{g_{0}}$ is the Laplace-Beltrami operator with respect to the metric
$g_0$ (with the sign chosen so that it is non-negative).
Thus equation (\ref{curvua}) can be written as:
\begin{equation}\label{eqU}
\Delta_{g_{0}}u=1-h^{2}e^{2u}-e^{-2u}\,\mbox{\rm det}_{g_{0}}(q).
\end{equation}
We now explain why (\ref{eqU}) has a unique smooth solution $u$ for each $h$ and $q$ fixed.
Since $h^2>0$ and $\mbox{\rm det}_{g_{0}}(q)\leq 0$ (recall that $q$ has trace zero)
we can find constants $u_{-}\leq u_{+}$ such that
\[1-h^{2}e^{2u_{-}}-e^{-2u_{-}}\,\mbox{\rm det}_{g_{0}}(q)\geq 0,\]
\[1-h^{2}e^{2u_{+}}-e^{-2u_{+}}\,\mbox{\rm det}_{g_{0}}(q)\leq 0.\]
It is well known that under these conditions (cf. \cite{KW1,KW2}) the semi-linear
elliptic equation (\ref{eqU}) admits a $C^{\infty}$ solution $u$ with
$u_{-}\leq u\leq u_{+}$ which must be unique by a maximum principle argument.

For each nonzero $h$, let $\mathcal S_h$ be the space of solutions of
(\ref{div2a}) and (\ref{curvua}). The discussion above shows that $\mathcal S_h$
is parametrized by pairs $(g_0,q)$ where $g_0$ is a metric of constant
curvature $-1$ and $q$ is the real part of a quadratic differential of the complex structure
determined by $g_0$. Of course, the diffeomorphism group $\mathcal D$ of $M$ and its identity component
$\mathcal D_{0}$ act on $\mathcal S_{h}$.

Ultimately we are interested in the flow $\phi$ up to smooth time changes and smooth conjugacy.
Having this in mind, we note that we only need to consider $\mathcal S_1$.
Indeed, if $(g,q,h)$ is a solution of (\ref{div2a}) and (\ref{curvua}) then
$(h^2\,g,h\,q,1)$ is also a solution and the corresponding thermostats are the same up to
a constant time change (cf. Lemma \ref{repa}).

Given $(g,q)\in \mathcal S_1$, we let $r^\pm=\pm 1+V(q)/2$. The smooth 2-dimensional bundles
spanned by $F$ and $H+r^{\pm}V$ give rise to two smooth transversal foliations $\mathcal F^{\pm}$.
The foliations are invariant under the thermostat $\phi$ defined by $(g,q)$.
By the results of Ghys \cite{Ghy2}, each $\mathcal F^{\pm}$ is smoothly conjugate to the weak foliation
of the geodesic flow of a metric $g_{\pm}$ of constant curvature $-1$.
Let $\mathcal T$ denote the Teichm\"uller space of the surface $M$ and note that each $g_{\pm}$ determines an element
$[g_{\pm}]\in \mathcal T$. Thus the map $\mathcal S_{1}\ni (g,q)\mapsto ([g_{+}],[g_{-}])\in \mathcal T\times\mathcal T$
induces a map
\[G:{\mathcal S}_{1}/{\mathcal D}_{0}\to \mathcal T\times\mathcal T\]
and since ${\mathcal S}_{1}/{\mathcal D}_{0}$ can be naturally identified with
$T(\mathcal T)$, the tangent bundle of $\mathcal T$, we get a map
\[G:T(\mathcal T)\to \mathcal T\times\mathcal T.\]
We do not really know that much about $G$. It clearly maps the zero section of $T(\mathcal T)$ onto the diagonal
of $\mathcal T\times\mathcal T$ and it is easy to see that
if $G([g_{0}],q)=([g_{+}],[g_{-}])$, then $G([g_{0}],-q)=([g_{-}],[g_{+}])$.

\medskip

\noindent{\bf Question.} Is $G$ a diffeomorphism?

\medskip

In the next subsection we will see that $G$ also takes values outside the diagonal of $\mathcal T\times\mathcal T$.

\subsection{Entropy production}

Consider the thermostat $\phi$ determined by an arbitrary function
$\la\in C^{\infty}(SM)$.
The next result is taken from \cite{DP3} and we include a proof for completeness sake.

\begin{Theorem} Let $p\in C^{\infty}(SM)$ be such that
$X(p)+HV(p)/k=0$ for some positive integer $k$,
and suppose
\[K-H(\la)+\la^{2}[(k+1)^{2}/(2k+1)]<0.\]
Then there exists $u\in C^{\infty}(SM)$ such that $F(u)=p$ if and
only if $p=0$.
\label{kten}
\end{Theorem}

\begin{proof}
We will use the following identity proved in \cite[Equation
(13)]{DP2}. Given $u\in C^{\infty}(SM)$ we have:

\begin{align}\label{id1}
2\int_{SM} Hu\, VFu\,d\mu
&=\int_{SM}(Fu)^2\,d\mu+\int_{SM}(Hu)^2\,d\mu\\
&-\int_{SM}(K-H(\lambda)+\lambda^2)(Vu)^2\,d\mu.\notag
\end{align}
Using that $X(p)+HV(p)/k=0$ and that $H$ and $X$ preserve the
Liouville measure we obtain:
\[\int_{SM}Hu\, V(p)\,d\mu=-\int_{SM}u\,HV(p)\,d\mu
=k\int_{SM}u\,X(p)\,d\mu=-k\int_{SM}X(u)\,p\,d\mu.\] Since
$X(u)=p-\la V(u)$ we derive
\[\int_{SM} Hu\,
VFu\,d\mu=-k\int_{SM}p^{2}\,d\mu+k\int_{SM}\la\,V(u)\,p\,d\mu.\]
Combining the last equality with (\ref{id1}) yields
\[(2k+1)\int_{SM}p^{2}\,d\mu-2k\int_{SM}\la V(u)\,p\,d\mu\]
\[+\int_{SM}(Hu)^2\,d\mu
-\int_{SM}(K-H(\lambda)+\lambda^2)(Vu)^2\,d\mu=0.\]
We may rewrite this equality as:
\begin{align*}
&\int_{SM}\left(\sqrt{2k+1}\,p-\frac{k\la\,V(u)}{\sqrt{2k+1}}\right)^2\,d\mu
\\&-\int_{SM}\left(K-H(\la)+\la^{2}\frac{(k+1)^{2}}{2k+1}\right)(V(u))^2\,d\mu
+\int_{SM}(Hu)^2\,d\mu=0.
\end{align*}
Combining this equality with the hypotheses we obtain $Hu=Vu=0$ which implies right away 
that $u$ must be constant.

\end{proof}

\begin{Remark}{\rm If $p(x,v)=q_{x}(v,\dots,v)$ where $q$ is a symmetric $k$-tensor, then
the condition $X(p)+HV(p)/k=0$ is just saying that $q$ has zero divergence (cf. proof of Lemma \ref{ldiv}).
For such a $p$ and $k=1$ it suffices to assume that $\phi$ is Anosov \cite{DP2}.
It is unknown if the Anosov hypothesis is enough for $k\geq 2$. The problem is open even for
geodesic flows. We refer to \cite{SU} for partial results in this direction when $k=2$.}
\end{Remark}

We now note that if $p\in H_{-k}\oplus H_{k}$, $k\geq 2$, then $X(p)+HV(p)/k=0$ is equivalent to saying that
$\eta_{-}p_{k}=0$, where $p=p_{-k}+p_{k}$ and
\[p+iV(p)/k=2p_{-k},\]
\[p-iV(p)/k=2p_{k}.\]
But the kernel of the elliptic operator $\eta_{-}$ in $C_{k}^{\infty}(SM)$ is a finite dimensional vector space which can be identified with the space
of holomorphic sections of the $k$-th power of the canonical line bundle. By the Riemann-Roch theorem, this space
has complex dimension $(2k-1)(\mbox{\rm genus}(M)-1)$ (for $k=2$ we get the holomorphic quadratic differentials).

Suppose $\phi$ is Anosov and let $\rho$ be the SRB measure.
{\it The entropy production} of the state $\rho$ is given by
\[e_{\phi}(\rho):=-\int \Div F\,d\rho=-\sum\,\mbox{\rm Lyapunov exponents}\]
where $\Div F$ is the divergence of $F$ with respect
to any volume form in $SM$. Recall that $V(\la)$ is the divergence of $F$ with respect to
$\Theta$.

D. Ruelle \cite{Ru0} has shown that $e_{\phi}(\rho)\geq 0$ 
and it is not hard to see (cf. \cite{DP2}) that $e_{\phi}(\rho)=0$ if and only if there is a smooth solution
$u$ to the cohomological equation $F(u)=V(\la)$.
Equivalently $\phi$ is volume preserving if and only if $e_{\phi}(\rho)=0$.
Theorem \ref{kten} gives a large class of thermostats with positive entropy production
just by taking $\la=\Re(p_{k})$, where $\eta_{-}p_{k}=0$. In this case $V(\la)=\Re (ikp_{k})$.
Note that $\phi$ is reversible if $k$ is odd so we get lots of new examples to which
the Fluctuation Theorem of G. Gallavotti and E.G.D. Cohen applies \cite{GC,GC1,Ga0,G}.

Let $(g,q)\in\mathcal S_1$ be a pair as in the previous subsection. Then $V(q)$ is the real part
of the restriction to $SM$ of a holomorphic quadratic differential. If $q$ is small enough then $\phi$ is an Anosov flow 
and by Theorem \ref{kten}, $\phi$ is volume preserving if and only if $q=0$.
This shows that the map $G$ above is interesting and that dissipative quasi-fuchsian flows do appear
as thermostats associated with quadratic differentials as claimed in the Introduction.

\end{document}